\documentclass{amsart}
\usepackage{amssymb}

\title[Moduli spaces of sheaves and Seiberg Witten invariants]
{Probing Moduli Spaces of Sheaves with Donaldson and Seiberg Witten Invariants}
\author {Rogier Brussee}
\address{ Telematica Instituut\\
          P.O. Box 589 \\
	    NL 7500 AN Enschede}
\email{rogier.brussee@telin.nl}
\dedicatory{In memory of Andrej Tyurin}

\makeatletter

\theoremstyle{plain}
\begingroup  
\newtheorem{Theorem}{Theorem} 
\newtheorem{Corollary}[Theorem]{Corollary} 

\newtheorem{Lemma}[Theorem]{Lemma} 

\endgroup

\theoremstyle{definition}

\theoremstyle{remark}
\newtheorem{Remark}[Theorem]{Remark}

\newtheorem{Acknowledgements}{Acknowledgments} 

\newcommand\url[1]{\hskip 0.2 ex plus 0.1 fil #1} 
\newcommand\txt[1]{\hbox{ #1 }}

\newdimen\arrowlen \arrowlen=2em
 
\newcommand\m@p[4]{\buildrel%
\hbox spread \arrowlen{\skip@=-0.33\arrowlen plus 1 fil\hskip\skip@$\m@th\scriptstyle #4$\hskip\skip@}%
\over%
{\mathord#1\mkern-6mu\cleaders\hbox{$\mkern-2mu\mathord#2\mkern-2mu$}%
\hfill\mkern-6mu\mathord#3}}
\newcommand\map{\m@p--\rightarrow} 
\newcommand\pam{\m@p\leftarrow--} 
\newcommand\equal{\m@p===} 
\newcommand\Map{\m@p==\Rightarrow} 
\newcommand\Pam{\m@p\Leftarrow==}
\newcommand\notdiv{\buildrel{\it/}\over{\smash|}}
\newcommand\rmmath[1]{\operatorname{#1}\nolimits}

\newcommand\Ext{\rmmath{Ext}}

\newcommand\sheafEnd{\operatorname{\mathcal{E}\mkern-3mu\mathit{nd}}\nolimits}

\newcommand\Pic{\rmmath{Pic}}

\newcommand\rank{\rmmath{rank}}

\newcommand\tr{\rmmath{tr}}
\newcommand\ord{\rmmath{ord}}
\newcommand\tensor{\otimes}
\newcommand\directsum{\mathop\oplus}
\newcommand\ch{\rmmath{ch}}
\newcommand\td{\rmmath{td}}
\newcommand\id{\rmmath{id}}

\let\slant=/
\newcommand\floor[1]{\left\lfloor #1 \right\rfloor}

\newcommand\iso{\cong}

\newcommand\R{{\mathbb R}}
\newcommand\C{{\mathbb C}}

\newcommand\Z{{\mathbb Z}}

\renewcommand\P{{\mathbb P}}
\renewcommand\O{{\mathcal O}}
\newcommand\F{{\mathcal F}}
\newcommand\numfrac[2]{\frac{#1}{#2}}
\newcommand\quart{{\numfrac14}}
\newcommand\half{{\numfrac12}}
\newcommand\threehalf{{\numfrac32}}

\let\next=\~\let\~=\tilde\let\tilde=\next
\let\next=\^\let\^=\hat\let\hat=\next

\def\<#1>{\left<#1\right>}
\let\trem@=\"\def\"#1{{\trem@{\if#1i\i\else#1\fi}}}

\newcommand\dual{{\scriptscriptstyle\vee}}
\newcommand\ddual{{\scriptscriptstyle\vee\vee}}
\renewcommand\({\left(}
\renewcommand\){\right)}

\newcommand\U{\mathrm{U}}
\newcommand\asd{{\mathrm{asd}}}

\renewcommand\L{{\mathcal L}}
\newcommand\M{{\mathcal M}}
\newcommand\Mbar{{\mkern4mu\overline{\mkern-4mu\M}}}

\newcommand\Mbartilde{{\widetilde\Mbar}}

\newcommand\Masd{\M_\asd}
\newif\ifcomment

\newcommand\Spin{{\mathrm{Spin}}}
\newcommand\Spinc{\Spin^c}

\newcommand\Kmin{K_{\min}}
\newcommand\Xmin{X_{\min}}
\newcommand\Pbar{{\bar\P}}

\renewcommand\L{{\mathcal L}}
\newcommand\FF{{\mathfrak F}}

\newcommand\odd{{\rmmath{odd}}}
\newcommand\km{{\rmmath{km}}}
\newcommand\sw{{\rmmath{sw}}}
\makeatother

\begin{document}
\begin{abstract}
We use Donaldson invariants of regular surfaces with $p_g >0$ to make quantitative 
statements about modulispaces of stable rank 2 sheaves. We give two examples: a quantitative existence theorem 
for stable bundles,  and a computation of the rank of the canonical holomorphic two forms on the moduli space. 
The results are in some sense dual to the Donaldson and O'Grady non vanishing theorems because they use 
the Donaldson series of the surface as {\em input}. Results in purely algebraic geometric terms 
can be obtained by using the explicit form of the Donaldson series of the surface. The Donaldson series are easy 
to compute using the Seiberg Witten invariants and the Witten conjecture which has recently been rigorously proved
by Feehan and Leness. 
\end{abstract}

\maketitle

\section{Introduction}

This paper uses gauge theory to get results about moduli spaces of stable sheaves on complex surfaces. There has 
long been an intimate relationship between gauge theory and the geometry of complex surfaces.  
Some of the earliest results about the Donaldson invariants concerned K\"ahler 
surfaces with $p_g >0$. By Donaldson's resolution of the Kobayashi Hitchin conjecture 
\cite[ch. 6]{DonaldsonKronheimer:FourManifoldsBook}, \cite{LuebkeTeleman:KobayashiHitchinBook}, one can 
compute their Donaldson invariants using complex geometry by determining the moduli space of stable 
rank 2 bundles and sheaves. Studying moduli spaces gave usefull qualitative 
non vanishing results for Donaldson invariants \cite{DonaldsonKronheimer:FourManifoldsBook}
\cite{OGrady:AnaloguesDonaldsonPolynomials} 
and provided explicit computations for special surfaces such as the elliptic ones 
(see \cite{FriedmanMorgan:FourManifoldsSurfacesBook} and references therein). 
Later, cut and paste techniques were developed to compute the invariants of a more general 
four-manifold \cite{MorganMrowkaRubermann}, \cite{FintushelStern:RationalBlowdown} that were also 
practical for e.g. the elliptic surfaces.  This work 
used computations for special complex surfaces like the elliptic surfaces and the K3 surface as 
as their input. At roughly the same time, Kronheimer and Mrowka organised all polynomials 
together in a generating series and proved a deep structure theorem for this Donaldson series
that revealed the essential content of the Donaldson polynomials (\cite{KronheimerMrowka:GaugeTheoryEmbeddedI} 
\cite{KronheimerMrowka:GaugeTheoryEmbeddedII}\cite{KronheimerMrowka:EmbeddedSurfaces}.
They, and through slightly different means Fintushel and Stern \cite{FintushelStern:DonaldsonSeries}, 
show that under the somewhat mysterious 
``simple type'' condition, the Donaldson series of a four manifold $X$ has the form of a universal 
power series in the intersection form and a finite number of ``basic classes'' 
$K_1, \ldots, K_n \in H^2(X,\Z)$ with $K_i \equiv w_2(X) \bmod 2$ (see theorem \ref{Theorem:DonaldsonSeries}). 

Soon after(and probably motivated by) Kronheimer and Mrowka's work, Seiberg Witten 
invariants were introduced by Witten \cite{Witten:MonopolesFourManifolds}. 
These also involved a finite number of characteristic classes in the sense above. 
The Seiberg Witten invariants are much more computable than the Donaldson Polynomials 
and they can easily be computed in full for complex surfaces 
\cite{Witten:MonopolesFourManifolds}\cite{FriedmanMorgan:AlgebraicSurfacesSeibergWitten}\cite
{Brussee:CanonicalClass}\cite{Nicolaescu:SeibergWittenBook}. Based on the physical heuristics, 
Witten conjectured  that the Kronheimer and Mrowka and Seiberg Witten 
basic classes coincide and the that the Donaldson series and the Seiberg Witten invariants 
mutually determine each other. Soon after Wittens work, Pidstrigatch and Tyurin proposed a mathematically 
rigorous approach to prove the Witten conjecture \cite{PidstrigatchTyurin:Cobordism} motivated by their previous 
work on the Spin polynomials \cite{PidstrigatchTyurin:SpinPolynomials}, \cite{Tyurin:DistinguishAlgebraicSurfaces}.
While conceptually very appealing, the program turned out to be very challenging technically.
In 2006, after an impressive series of papers spanning almost 10 years, 
Feehan and Leness finally cracked the problem in the generality needed here 
\cite{FeehanLeness:PU2MonopolesI}\cite{FeehanLeness:PU2MonopolesII}\cite{FeehanLeness:PU2MonopolesIII}\cite
{FeehanLeness:PU2MonopolesSurvey}\cite{FeehanLeness:GeneralFormula}\cite{FeehanLeness:WittenConjecture}. 
More detailed expositions of the physical arguments for the Witten conjecture have also become available 
(see e.g. \cite{MooreWitten:IntegrationUPlane}\cite{LabastidaMarino:TopologicalQFTFourManifoldsBook}\cite
{Vajiac:WittenConjecture}).

The net result of these developments is that the Donaldson polynomials have become very
computable invariants and thus can be used to get {\em a priori} information about moduli spaces of stable sheaves 
i.e. with information flowing in the opposite direction compared to the  early days of Donaldson theory! 
In this paper we consider two examples of this strategy.  

\subsection{Algebraic geometric results for moduli spaces of sheaves}

Let $X$ be a regular complex algebraic surface of general 
type with $p_g >0$ i.e. with $b_+ \ge 3$ and $b_1 = 0$. Fix a complex line bundle $\L \in \Pic(X)$ 
with $c_1(\L) = L$, and a polarization $H$. For $k \in \Z$, let $d(L,k) = 4k - L^2 -3(1+p_g)$ be the 
complex virtual dimension of the moduli space of stable bundles with $c_1 = L$, $c_2 = k$. Finally let 
$L_{\min}$ be the component of $L$ in the cohomology of the minimal model, and $\odd(L)$ the number 
of $(-1)$-curves $E$ with $L\cdot E$ odd.

First we prove that semi stable bundles with small virtual dimension exist. The following theorem
 is a direct consequence of theorem \ref{Theorem:ExistenceStableBundles} 
and the Donaldson series of surfaces of general type in theorem \ref{Corollary:DonaldsonSeriesSurfaces}

\begin{Theorem} \label{theorem:ExistenceGeneraType}
For a surface of general type as above there exists an  $H$ semi stable bundle $V$ 
with $\det(V) = \L$ and $c_2(V) = k$ where $k$ satisfies  
$$
	d(L,k) = 4k - L^2 -3(1+p_g) \le \odd(L) + 3
%
$$
\end{Theorem}

This result seems to become less sharp as we blow up the surface more often. However,
if the the polarization $H\cdot L$ is odd and is sufficiently close to a class pulled 
back from the minimal model (i.e. $H\cdot E_i < \epsilon (H^2)^{1/2}$ for all exceptional curves and 
sufficiently small $\epsilon$), then 
$$
	d(L,k) \ge \odd(L) - 3(1+p_g).
$$
This follows from the Bogomolov inequality and the fact that for every stable bundle on the blow up, 
the double dual of the push forward of the minimal model is a semi-stable bundle with 
at most the same second Chern class \cite{Brussee:BlownupSurfaces}.   

Likewise we have a result for elliptic surfaces that follows from theorem \ref{Theorem:ExistenceStableBundles} 
and the shape of the Donaldson series of elliptic surfaces in theorem \ref{Corollary:DonaldsonSeriesSurfaces}. 
We only state a result for the minimal case.

\begin{Theorem}
 Let $X \to \P^1$ be a regular elliptic surface with $p_g >0$, 
general fibre $F$ and multiple fibres $F_{p_1} \ldots F_{p_n}$ with 
multiplicities $p_1, \ldots p_n$. Assume that $2 \notdiv \gcd(p_1, \ldots p_n)$.
Let $L$ be a vertical divisor (i.e. a rational multiple of $F$) and choose   
a polarization $H$. Then there is $H$-semi-stable bundle $V$ 
with $c_1(V) = L$ and $c_2(V) = k$ with
$$
	d(L,k) \le  n + p_g - 1
$$ 
\end{Theorem}

As a second application we compute the generic rank of the canonical two forms on the moduli space 
of Gieseker semi stable sheaves (see section \ref{Section:GenericRankTauForm}) defined by Mukai and Tyurin 
\cite{Mukai:SymplecticStructureModuliSheavesK3Abelian}\cite
{Tyurin:SymplecticStructureModuliSpaceVectorbundlesAlgebraicSurfaces}.   
Theorem \ref{Theorem:GenericRankTauFormGTElliptic} extends 'O Grady's non vanishing result 
\cite[theorem 2.4]{OGrady:AnaloguesDonaldsonPolynomials}.
As far as I know, the best existing result is a a refinement by Jun-Li \cite{JunLi} who proves 
that if $X$ is regular of general type with $p_g >0$, and the linear system $|K|$ contains a reduced 
irreducible divisor $D$, then the two form on the the moduli space generically has the maximal possible 
rank $\floor{\half d(L,k) }$. O'Grady's and Li's computation is entirely algebraic geometric. O'Grady showed that 
this implies the non vanishing of the Donaldson polynomials $q_{L,k}$ evaluated on $d(L,k)$ copies 
of $\omega + \bar\omega$ for even dimensional moduli space and on $d(L,k) - 1 $ copies of 
$\omega + \bar \omega$ and one copy of the hyperplane section $H$ in the odd case. 
Note that $D$ is necessarily reducible (even disconnected) if $X$ is non minimal or not 
a surface of general type. Using the implication in the other direction and the Donaldson series of 
a surface of general type, the theorem below  gives us a result that
removes {\em a posteriori} all conditions on the canonical divisor in the minimal case.

It seems natural that a purely algebraic proof of \ref{Theorem:GenericRankTauFormGTElliptic} below is possible, 
especially in the minimal general type case, but that something like numerical connectivity of the 
canonical divisor \cite[prop. 6.2]{BarthPetersVanDeVen:ComplexSurfacesBook} must be used. 
However such a proof seems to be unavailable and it would have been very useful 
in applications to 4-manifold theory before Seiberg Witten theorie was developed 
\cite{Kronheimer:ObstructionIntersectionPoints}
\cite[prop. 9.6]{KronheimerMrowka:EmbeddedSurfaces}, \cite{Brussee:(-1)-SphereConjecture},
\cite{Brussee:RemarksKronheimerMrowkaClasses}. It introduced existence of an irreducible canonical divisor as a 
condition in several statements which naturally raised the question whether this condition 
was just a technical assumption reflecting our lack of technical knowhow
or pointed at a deep differentiable obstruction to the existence of such a divisor. 
Apparently both is true, but only the obvious reasons for reducibility, 
non minimality and non general type, play a role. 
Theorem \ref{Theorem:GenericRankTauFormGTElliptic} below in the minimal general type case is referred to in 
\cite[above Theorem 10.4.4]{HuybrechtsLehn:SheavesBook}.
It is a direct consequence of the shape of the Donaldson series (theorem \ref{Corollary:DonaldsonSeriesSurfaces})
and theorem \ref{Theorem:GenericRankTauForm}.
  
\begin{Theorem}\label{Theorem:GenericRankTauFormGTElliptic}
For an algebraic surface $X$ as above with polarisation $H$, let $\sigma \in H^0(X,K_X)$ be a holomorphic two form,
and let $\tau(X)$ be the associated two form on the stable locus of the moduli space of Gieseker stable sheaves 
$\Mbar_H(L,k)^s$.  
\item Assume that $X$ is a surface of general type, then for $k \gg 0$,  
$$
	\rank\(\tau(\omega)\) = \floor{\half\(d(L,k) - \odd(L)\)}
$$
\item Assume that $X  \to \P^1$ is a regular elliptic algebraic surface 
with $p_g >0$, general fibre $F$ and multiple fibres $F_{p_1} \ldots F_{p_n}$ with 
multiplicities $p_1, \ldots p_n$. Also assume that $2 \notdiv \gcd(p_1, \ldots p_n)$.
Let $L$ be a vertical divisor (i.e. a rational multiple of $F$). Then for $k \gg 0$, 
$$
	\rank\(\tau(\omega)\) = \floor{\half\(d(L,k) - (p_g -1) \)} = 2k - 2p_g - 1
$$
\end{Theorem}

\subsection{Organisation of the paper}
The organization of this paper is as follows. In section \ref{Section:Preliminaries} 
we first do some preliminaries. We recall the Donaldson series,
state the Witten conjecture, write down the Donaldson series of for surfaces of different Kodaira dimension, 
and explain the construction of the holomorphic form on the moduli space. In section \ref{Section:ExistenceStableBundles} 
we then prove the existence of semi stable bundles with low virtual dimension. 
The only problem here is properly dealing with semi stable bundles which we handle by blowing up the surface. 
In the last section \ref{Section:GenericRankTauForm} we compute the 
rank of the two form on the moduli space given its Donaldson series.
Semi stable sheaves are again handled by blowing up the surface but the computations are less straight forward than
in section \ref{Section:ExistenceStableBundles}.

\begin{Acknowledgements}
This work was started while I worked at Bielefeld University and the SFB 343. 
Countless talks with Stefan Bauer, Manfred Lehn and Victor Pidstrigatch 
on gauge theory and algebraic geometry have shaped my ideas on the subject. The book by Manfred Lehn 
and Daniel Huybrechts which I saw growing there, has been particularly 
useful for this paper. Talks with Vicente Mu\tilde noz and Matilde Marcolli have also been influential. 
Thanks also to Mettina Veenstra and the Telematica Instituut for providing the opportunity to finish this paper.
Finally, special thanks to Carolien Kooiman for providing mental support.

This paper is dedicated to Andrej Tyurin. It was a shock to read that he passed away. I have been 
heavily influenced by his papers on gauge theory and complex geometry. He and Victor were also much closer on track
then anybody (or at least I) suspected with their Spin polynomial invariants,  
and it is tribute to their insight that their program to prove the Witten conjecture has finally been realized. 
Most of all though, I fondly remember Andrej as a man who remained young at heart. It is a privilege to have known him. 
\end{Acknowledgements}

\section{Preliminaries and Notation}\label{Section:Preliminaries}

This section introduces notation and collects some results. 
The material is either well known or just slightly massaged in 
a form appropriate for use in the last two sections.

\subsection{The Donaldson series}\label{SubSection:DonaldsonSeries} 
Mainly to fix notation we first recall Kronheimer and Mrowka's  work on the 
existence and shape of the Donaldson series \cite{KronheimerMrowka:EmbeddedSurfaces}.

Let $X$ be a closed 4-manifold  with odd $b_+ >1$, $b_1 = 0$. Fix a line bundle 
$\L$ with $c_1(\L) = L$ and fixed connection. Consider the moduli space $\Masd(\L, k)$, 
of irreducible ASD connections on a $\U(2)$ 
bundle $V$ with $c_1(V) = L$, $c_2(V) = k$ and determinant connection equal to the 
fixed connection on the line bundle $\L$. The moduli space has virtual real dimension 
$2d$ where for $b_1 = 0$, $d = d(k,L)$ is given by
$$
	d = 4k - L^2 - \threehalf(b_+ + 1).
$$
The Donaldson polynomial $q_{L,k}$ is constructed from  (a quarter of) the first Pontrijagin
class of the universal $PU(2)= SO(3)$ bundle over $\Masd(\L,k) \times X$ using the 
the slant map construction \cite[def. 5.1.11]{DonaldsonKronheimer:FourManifoldsBook}. 
It is considered as a homogeneous polynomial of degree $d$ on 
$H_2(X) \oplus H_0(X)$, by giving $\Sigma \in H_2(X)$
degree 1, and the homology class of a point $x$ degree 2. 
Using this convention, a manifold is of Donaldson simple type iff its Donaldson polynomials
satisfy
$$
	q_{L,k}(\Sigma^{d-4},x^2) = 4 q_{L,k-1}(\Sigma^{d-4})  
                           \qquad d = 4k -L^2 -\numfrac32 (1+ b_+).
$$
This condition seems to depend on the cohomology class $L$ but it is known that the simple 
type condition is satisfied for all classes $L$ iff it is satisfied for just one 
\cite[lemma 7.37]{KronheimerMrowka:EmbeddedSurfaces}. 
The simple type condition can also be formulated in terms of the relation between 
the polynomials of and the manifold  and those of the ``blow up'', the connected sum $ X\#\Pbar$ 
of the manifold with $\C\P^2$ with the opposite orientation.
For simple type manifolds the Donaldson series $q_L$ is then defined as an element of 
$(S^\bullet H_2(X))^\dual$ by
\begin{equation}
\label{Eq:DonaldsonSeriesDefinition}	
			q_L(X) = \sum_{d = 0}^\infty \numfrac{1}{d!}\,q_{d,L}.
\end{equation}
where
$$
	q_{d,L} = \begin{cases} 
			q_{L,k}                     & \txt{if} d = 4k-L^2 - \threehalf (1+ b_+), 
\\
			\half q_{L,k}(x,-)          & \txt{if} d = 4k -2 -L^2 - \threehalf (1+ b_+),
\\
			0                           & \txt{otherwise}. 
		    \end{cases}
$$
In particular the series is even or odd depending on the parity of $L^2 + (b_+ + 1)/2$. 

\begin{Theorem} \label{Theorem:DonaldsonSeries}(Kronheimer Mrowka)
Let $Q$ be the intersectionform on $X$. there exist a finite number of {\em Kronheimer-Mrowka basic classes}  
$K_1, \ldots, K_n \in H^2(X,\Z)$ such that
\begin{equation}
\label{Eq:DonaldsonSeries}		
		q_L(X) = e^{Q/2}\sum_{i=1}^n  (-1)^{(L^2 + K\cdot L)/2} \km(X,K_i) e^{K_i}.
\end{equation}
The basic classes are characteristic i.e. reduce to $w_2(X) \bmod 2$. 
\end{Theorem}

\begin{proof}
\cite[Theorem 1.7]{KronheimerMrowka:EmbeddedSurfaces} 
\end{proof}

The rational numbers $\km(X,K_i)$ will be called the Donaldson multiplicity of the 
Kronheimer Mrowka basic class $K_i$. Since the  Donaldson series are even or odd 
$\km(X,-K_i) = \pm (X,K_i)$, and the basic classes $K_i$ come in pairs 
differing by a sign.

\subsection{Seiberg Witten invariants and the Witten Conjecture}\label{SubSection:SeibergWittenAndWittenConjecture}

On the Seiberg Witten side,  let $S$ be a $\Spinc$ structure of a four manifold $X$.
For every such 
$\Spinc$ structure the Seiberg Witten invariant $\sw(X,S)$ of the pair $(X,S)$ 
is defined in terms of the solutions of the monopole equations for spinors in $S$ 
\cite{Witten:MonopolesFourManifolds},\cite{Brussee:CanonicalClass},\cite{Nicolaescu:SeibergWittenBook}. 
The Seiberg Witten multiplicity of the pair $(X,K)$ where $K \in H^2(X,\Z)$ is a cohomology class,
is then defined by
$$
	\sw(X,K) = \sum_{c_1(S^+) = K} \sw(X,S).
$$
The finitely many classes $K_i$ with $\sw(X,K_i) \ne 0$ are called the Seiberg Witten basic classes. 
It is easy to prove that $\sw(K) = \pm \sw(-K)$ 
so $K$ is a Seiberg Witten basic class, iff $-K$ is. 
Since $c_1(S^+) \equiv w_2(X) \bmod 2$ for every $\Spinc$-structure $S$, 
the Seiberg Witten basic class $K$ is characteristic by construction. 

A four-manifold is of {\em Seiberg Witten simple type} if all basic classes $K$ satisfy 
$$
	K^2 = (2e + 3\sigma)(X),
$$
where $e$ and $\sigma$ are the topological Euler characteristic and the signature
of $X$ respectively. Equivalently, a four manifold is of simple type, iff only 
moduli spaces of solutions to the monopole equation with virtual dimension 0 give rise to 
non trivial Seiberg Witten invariants.

We now state the Witten conjecture \cite{Witten:MonopolesFourManifolds} which has recently been proven 
by Feehan and Leness.  

\begin{Theorem}\label{Theorem:WittenConjecture} (Feehan Leness, Witten conjecture) 
Let $X$ be a 4 manifold  with odd $b_+ \ge 3$ and $b_1 = 0$. 
\begin{enumerate}
\item If $X$ is of Seiberg Witten simple type then it is of Kronheimer Mrowka simple type.
\item the Seiberg Witten basic classes and the Kronheimer Mrowka basic classes coincide. More precisely
\begin{equation}
\label{Eq:SW=KM}	\km(X, K) =  2^{2 + \quart(7e + 11\sigma)(X)} \sw(X,K).
\end{equation}
\end{enumerate} 
\end{Theorem}

\begin{proof}
\cite[theorem 1.2]{FeehanLeness:WittenConjecture}. The proof is a clever combination of Feehan
 and Leness's earlier results on the
Pidstrigatch Tyurin $PU(2)$-cobordism program and fixing universal functions by comparing their results to a 
large supply of examples with only one basic class (in both senses) constructed by Fintushel and Stern.
\end{proof}

\subsection{Seiberg Witten invariants and Donaldson Series of Surfaces}
\label{SubSection:SeibergWittenDonaldsonSeriesSurfaces}

What we need here on the Seiberg Witten invariants of K\"ahler surfaces  is contained 
in the following theorem.
 
\begin{Theorem} \label{Theorem:BasicClassesSurfaces}
Let $X$ be a K\"ahler surface with $p_g >0$ with canonical class $K_X \in H^2(X,\Z)$. 
\begin{enumerate}
\item The surface $X$ is Seiberg Witten simple type.
\item If $X$ is a surface of general type, then the basic classes are all classes of the form 
$$
	\pm \Kmin + \sum_i \pm E_i
$$ 
where $\Kmin$ is the canonical class of the minimal model and the $E_i$ run over all
$(-1)$-curves. 
Moreover $\km(X,-K_X) = 1$.
\item
If $X \to C$ is a minimal elliptic surface with $p_g> 0$, general fibre $F$ and 
multiple fibres $F_1 \ldots F_n$ of multiplicity $p_1, \ldots p_n$,
then the basic classes are of the form 
$$
	K = -K_X + 2 (d F + \sum a_i F_i)
$$
where $0\le d \le (p_g-1) + g(C)$ and $0 \le a_i < p$. Moreover
\begin{equation}
\label{Eq:SWElliptic}	\sw(K) = \sum_{K  = -K_X + 2(dF + \sum a_i F_i)} (-1)^{d} \binom{p_g-1 + g(C)}{d}
\end{equation}
where the sum runs over all possible choices for $\half(K +K_X)$ in $H^2(X,\Z)$.   
\end{enumerate}
\end{Theorem}

\begin{proof}
The general type case see \cite[Cor. 31, lemma 11, prop. 41]{Brussee:CanonicalClass} or 
\cite[prop.2.5]{FriedmanMorgan:AlgebraicSurfacesSeibergWitten}. It is in essence  already in Wittens paper
\cite{Witten:MonopolesFourManifolds}. See also \cite[prop. 3.3.1]{Nicolaescu:SeibergWittenBook}.
For the elliptic case see \cite[Prop. 42] {Brussee:CanonicalClass} or\cite[Prop. 4.4]{FriedmanMorgan:ObstructionBundles}
See also \cite[ch. 3.3.20]{Nicolaescu:SeibergWittenBook} for the simply connected case.
\end{proof}

We now introduce the following ad hoc notation. For a pair of a K\"ahler surface $X$ and a class 
$L \in H^2(X,\Z)$ we let $\odd(L)$  be the number of $(-1)$-curves $E_i$, with $E_i\cdot \L$ odd. 
W.l.o.g. we number the $(-1)$-curves in such a way, that the first $\odd(L)$ 
$(-1)$-curves have odd intersection with $L$. 

\begin{Corollary} \label{Corollary:DonaldsonSeriesSurfaces}
Let $X$ be regular complex surface with $p_g >0$ (i.e $b_+ \ge 3$ odd and $b_1 = 0$), 
and $L \in H^2(X,\Z)$.
\begin{enumerate}
\item 
If $X$ is of general type, then the Donaldson series is given by
$$
          q_L = q_0e^{Q/2}(e^{-\Kmin} + (-1)^{1+ p_g + L_{\min}^2} e^{\Kmin})
                               \prod_{i\le \odd(L)}\sinh(E_i)\prod_{i>\odd(L)}\cosh(E_i)
$$
where the constant is given by
$$
	q_0 = (-1)^{\half(L_{\min}^2 - K_{\min}\cdot L_{\min})} 2^{2 + \quart(7e + 11\sigma)(\Xmin)}.
$$
\item
If $X \to \P^1$ is a minimal elliptic surface with general fibre $F$ and
multiple fibres $F_1 \ldots F_n$ of multiplicity $p_1, \ldots p_n$ such that $2 \notdiv \gcd(p_1, \ldots p_n)$ 
and and $L$ is a rational multiple of $F$ then
$$
	q_L= q_0 e^{Q/2} \frac{\sinh^{p_g-1 + n}(-F)}{\prod_i \sinh(-F_i)}
$$
where the constant is given by
$$
	  2^{2 + \quart(7e + 11\sigma)(X) + (p_g -1)}.
$$
\end{enumerate}
\end{Corollary}

\begin{proof}
For the general type case this is an immediate consequence of the Witten conjecture (theorem 
\ref{Theorem:WittenConjecture}) 
and theorem \ref{Theorem:BasicClassesSurfaces}. 
However, it is easier to start with the minimal case and use the blow up formulas for the Donaldson polynomials
(\cite{FintushelStern:BlowupFormulas} or \cite[lemma 9.2, theorem 7.23]{KronheimerMrowka:EmbeddedSurfaces} for 
the simple type case). 
For minimal surfaces of general type, the corollary follows immediately, because in this case $K = -\Kmin$ 
is the only Seiberg Witten basic class up to sign and it has multiplicity $1$. 
Therefore, by the Witten conjecture, it is the only Kronheimer Mrowka basic class, 
and so the Donaldson series has the form given by \eqref{Eq:DonaldsonSeries}. For the general case, we decompose $L$ as 
$$
	L = L_{\min} + \sum_{i=1}^{\odd(L)} E_i
$$
where $L_{\min}$ is pulled back from the minimal model, and use the blow up formulas for the Donaldson Polynomials.

The Donaldson series for simply connected elliptic surfaces with two multiple fibres has been 
computed by Fintushel and Stern  without recourse to the Witten conjecture,
see \cite[theorem 1.1]{FintushelStern:RationalBlowdown}.
For the more general regular elliptic case, we have made simplifying assumptions. 
Since $L$ is vertical $\half(L^2 + L\cdot K_X) \equiv 0 \bmod 2 $ and we don't worry 
about signs in the Donaldson series which might change $\sinh$ in $\cosh$. 
The assumption on the $\gcd$ further implies that there is no two torsion in the cohomology so the sum
in equation \eqref{Eq:SWElliptic} collapses to a single summand. 

Now writing $K_X = (p_g -1)F + \sum_i (p_i -1) F_i$, we see that by the Witten conjecture \ref{Theorem:WittenConjecture} 
and by the explicit Seiberg Witten multiplicities in theorem \ref{Theorem:BasicClassesSurfaces}, the Donaldson series is 
\begin{align*}
	q_L &=  q'_0 e^{Q/2}\sum_{0 \le d \le (p_g -1), 0\le a_i \le (p-1)} 
                                   (-1)^d \binom{p_g -1}{d}e^{(-(p_g -1) + 2d)F + \sum_i(-(p-1) + 2a_i)F_i} 
\\
	    &=q'_0 e^{Q/2} (e^{-F} - e^F)^{p_g-1} \prod_i (\sum_{-(p-1) \le a_i \le (p-1)} e^{a_iF_i})
\\
	    &=q'_0 e^{Q/2} 2^{p_g - 1} \sinh^{p_g-1}(-F)\prod_i \sinh(-pF_p)/\sinh(-F_p)
\end{align*}
where $q'_0 = 2^{2 + \quart(7e + 11\sigma)(X)}$. Since $pF_p = F \in H^2(X,\Z)$ we are done.
\end{proof}

\subsection{Holomorphic Symplectic Forms on the Moduli Space of Sheaves}
\label{SubSection:TauForm}

A holomorphic two form $\omega \in H^0(K)$ on an algebraic surface $X$ defines a  holomorphic two form $\tau(\omega)$ 
on the moduli space of stable sheaves \cite{Mukai:SymplecticStructureModuliSheavesK3Abelian}\cite
{Tyurin:SymplecticStructureModuliSpaceVectorbundlesAlgebraicSurfaces}. 
We follow \cite[section 10.3]{HuybrechtsLehn:SheavesBook}. In a point $[\F]$ corresponding to a 
stable sheaf in the moduli space, the form is equal to  
$$
	\tau(\omega)_{\F} = \int_X \ch_2(\F)\wedge \omega, 
$$
where $\ch_2(\F)$ is half the trace of the square of the Atiyah class
$$
	A(\F) \in Ext_X^1(F,  \Omega^1_{X \times \Mbar}|_{X\times[\F]} \tensor \F ) ),
$$
\cite[section 10.1.5]{HuybrechtsLehn:SheavesBook} which lives in 
$$
	H^2(X,\Omega^2_{X \times \Mbar}|_{X\times [\F]})  = 
                   \directsum_{i=0}^2 H^2(X, \Omega_X^{2-i}) \tensor \Omega^i_{\Mbar}|_{[\F]}.
$$ 
The operation $\int_X \omega \wedge -$ is a convenient and suggestive notation 
for the cohomological operation of taking cup product with $\omega \in H^0(X,K_X)$, which 
for dimension reasons kills all but the $i=2$ component in the above direct sum, and 
using the canonical isomorphism $H^2(X,K_X) \iso \C$ to to get a class in 
$$
	H^2(X, K_X) \tensor \Omega^2_\Mbar|_{[\F]} \iso  \Omega^2_{\Mbar}|_{[\F]}.
$$ 

If we have a universal family $\FF$ of sheaves over $\Mbar^s$, the locus of the moduli space
corresponding to the Gieseker stable sheaves, we can easily globalize the above construction to
$$
	\tau(\omega) = \int_X \ch_2(\FF)\wedge \omega. 
$$
which is a holomorphic two form.  Such a universal family may not exist for all polarizations $H$ 
but there always exists a quasi universal sheaf 
$\FF$. A sheaf $\FF$ over 
$X \times \Mbar^s(\L,k)$ is quasi universal if for all connected flat families of Gieseker stable sheaves 
$\FF_S$ with $c_2 =k$ and determinant $\L$, parametrized by a connected scheme
$S$ there is a classifying map $f: S\to \Mbar^s(\L,k)$ and a vector bundle 
$W$ on $S$, such that $\F_S \tensor W \iso f^* \FF$. 
Clearly the rank of $W$ depends only on the connected component of $\Mbar^s$.
%
The quasi universal family still defines a two form on $\Mbar^s$ by
$$
	\tau(\omega) = \frac1{\rank(W)} \int \ch_2(\FF) \wedge \omega 
$$
independent of the choice of the quasi universal family. 
\cite[lemma 10.3.4 and Proposition 10.4.1]{HuybrechtsLehn:SheavesBook}. 

Mutatis mutandis the same construction 
can be made for antiholomorphic two forms to get a antiholomorphic 
two form on the moduli space. Because $\ch_2(\F)$ is a real class, it is easy to see that 
$$\label{conjugation}
	\overline{\tau(\omega)} = \tau(\bar\omega).
$$

\begin{Remark}
In Dolbeault cohomology, the curvature of a bundle with 
hermitian connection represents the Atiyah class up to the usual factor $2\pi i$ 
(see \cite[lemma 3.1, theorem 3.2]{OGrady:AnaloguesDonaldsonPolynomials}). Thus if we choose a 
connection $\nabla$ on a universal family 
of bundles over the product of (an open neighborhood in) the moduli space of stable bundles and the surface $X$,   
then in Dolbeault cohomology, the construction above is 
$$
	\int_{X}(\tr F^2(\nabla)) \wedge \omega
$$
where the integral over $X$ is fibre wise and the result is a 
two form on the moduli space. With the modification of dividing by $\rank(W)$ as above
we can also use this differential geometric construction with a quasi universal family.  

We also recall the interpretation of the rank of this two form which is the starting point for the algebraic
geometric approach to the O'Grady non vanishing theorem, from 
\cite[proof theorem 2.4, 2.7]{OGrady:AnaloguesDonaldsonPolynomials}\cite[Proposition 10.4.1]{HuybrechtsLehn:SheavesBook}. 
Consider a surface $X$ with $p_g > 0$ and $\omega \in H^0(K)$ a holomorphic two form 
vanishing on a divisor $D$. Let $V$ be a stable vector bundle on $X$ with $H^2(\sheafEnd_0(V)) = 0$. 
Then the following are equivalent:

\begin{enumerate}
\item the map induced by multiplication with $\omega$ 
$$	
	m:H^1(\sheafEnd_0(V)) \map{\cdot\omega} H^1(\sheafEnd_0(V)\tensor K)
$$
has a kernel and cokernel of dimension $e$.
\item 
$$
      h^0(\sheafEnd_0(V)\tensor K|_D) = e.
$$
\item
The form $\tau(\omega)^n \ne 0$ for $n \le \floor{d_{L,k} - e)/2}$ at the point $[V]$ in the moduli space 
\end{enumerate} 
\end{Remark}

\section{Existence of stable rank 2 bundles}\label{Section:ExistenceStableBundles}

We will first give an existence proof for semi stable bundles. 
The two main problems with using the Donaldson polynomials directly, is that in general the moduli spaces have
a geometric dimension larger than the virtual one and that we have to be careful about handling semi stable bundles.
In addition, we have to be careful with proper sheaves 
(as opposed to vector bundles) as they give rise to points on the boundary of the Uhlenbeck compactification 
of the ASD moduli space. We avoid these problems by going to the limit of small $c_2$ , and considering 
the blow up of the surface. 

\begin{Theorem}\label{Theorem:ExistenceStableBundles}
Let $(X, \Omega)$ be a complex k\"ahler surface with $p_g >0$, and $b_1 = 0$, and such that 
the Donaldson series $q_L$ is of order $n$ i.e. has the form 
$$
	q_L =  \sum_{d \ge n} \frac{1}{d!} \, q_{d,L} 
$$
where $q_{n,L} \ne 0$.  Then there is a $\Omega$-semi stable 
bundle with
$$
	d(L,k) = 4k - L^2 -3(1+p_g) \le n+2
$$	
\end{Theorem}

\begin{proof}
Since the Donaldson series is even or odd, it is trivially true that 
$-L^2 -3(1+p_g) \equiv n \bmod 2$. The inequality of the theorem is therefore trivially equivalent to
$$
	d(L,k)  \le \begin{cases}
	                  n      & \txt{if} -L^2 -3(1+p_g) \equiv n \bmod 4
\\
			      n+2    & \txt{if} -L^2 -3(1+p_g) \equiv n+2 \bmod 4
		       \end{cases}
$$   
We prove the case $-L^2 -3(1+p_g) \equiv n \bmod 4$ in detail, and indicate what has to be 
changed in the other case.

We first consider the case of K\"ahler forms $\Omega$ such that for all  $F \in \Pic(X)$ with 
$-(L- 2F)^2 - 3(1+ p_g) \le n$, we have
\begin{equation}
\label{Eq:NotOnWall}		\Omega\cdot (L - 2F) \ne  0,
\end{equation}
i.e. we assume that $\Omega$-semi stability is the same as $\Omega$-stability for all sheaves 
in the Gieseker compactification 
of the moduli space. In differential geometric terms this means that there does not exist a reducible ASD 
connections (or equivalently Hermite Einstein metrics) in the Uhlenbeck compactification. 
The "not on a wall condition'' \eqref{Eq:NotOnWall} is trivially 
satisfied if $\Omega\cdot L \notin 2\Omega\cdot \Pic(X)$ e.g. 
if $\Omega$ is an integral class and $\Omega\cdot L$ is odd.

Now suppose that there are no stable bundles with $d(L,k) \le  n = \ord(q_L)$. 
Then by the resolution of the Kobayashi Hitchin conjecture, there are no irreducible Hermite Einstein metrics, 
or equivalently, irreducible ASD connections on vector bundles $V$ with $d(L,k) \le n$ and 
fixed determinant connection on the differentiable line bundle underlying  $\L$. Therefore the Hodge metric is a 
perfectly good ``generic'' metric to compute the Donaldson invariants $q_{L,k}$ up to and including 
virtual dimension $n$, because all the moduli spaces in the Uhlenbeck compactification, being empty, 
have the expected dimension. 
Moreover, obviously
$$
	q_{L,k} = 0.
$$
However if we choose $k$ such that $d(L,k) = n$, which we can 
because $-L^2 -3(1+p_g) \equiv n \bmod 4$, then   
$$
	q_{L,k}|_{H_2(X)} =  n! \cdot \txt{degree $n$ part of} q_L \ne 0,   
$$
a contradiction. 

Essentially the same argument goes through in the other case $-L^2 -3(1+p_g) \equiv n +2 \bmod 4$, 
except we have to choose $k$ such that $d(L,k) = n+2$ and we get a contradiction from 
$$
	q_{L,k}(x, -)|_{H_2(X)} \ne 0.
$$

To get rid of the condition \eqref{Eq:NotOnWall}, we first consider the blow up of the surface $\sigma:\~X \to X$ 
and choose $\~L= L + E$, where $E$ is the newly introduced exceptional curve. 
We also choose a K\"ahler metric with K\"ahler form  $\~\Omega= \sigma^*\Omega - \epsilon E$ 
\cite[p. 186]{GriffithHarris:AlgebraicGeometryBook} where $\epsilon$ is small and such that 
$$
	\~\Omega\cdot \~L = \Omega\cdot L + \epsilon \notin 2 \~\Omega\cdot \Pic(\~X)
$$
e.g. if $\Omega$ is integral and $\epsilon =1/2N$ for $N \gg 0$. 
Note that Donaldson invariants are {\em defined} by counting ASD connections on the connected sum with $\Pbar$ to avoid 
the problem of reducible connections. Thus the blow up procedure (differentiably taking the connected sum with 
$\Pbar$) to handle reducible semi stable bundles is quite natural from our present point of view. 
 
By the blowup formula for the Donalson series 
(\cite{FintushelStern:BlowupFormulas} or \cite[lemma 9.2, theorem 7.23]{KronheimerMrowka:EmbeddedSurfaces})
$$
	q_{\~L} = \sinh(E)e^{-EE}q_L = Eq_{n,L}/n! + O(n+3),
$$
i.e. $q_{\~L}$ is of order $n+1$. Thus, by the above we find a stable bundle $\~V$ on $\~X$ with $c_1(\~V) = \~L$ and
$$
	 d(\~L,k) = d(L,k) + 1 \le (n + 1) + 2.
$$
The pushforward along the blowup $\F = \sigma_*(\~V)$ has $c_1(\F) = L$, 
and $c_2(\F) = k$ (c.f. the computations in 
\cite{Brussee:BlownupSurfaces}). After taking the double dual if necessary, which can only decrease $c_2$, 
we get a semi stable bundle $V = \F^\ddual$ with $d(L,k) \le \odd(L)$. 
\end{proof}

\section{Generic rank of two forms on the moduli space}\label{Section:GenericRankTauForm}

We now compute the generic rank of the canonical holomorphic two form $\tau(\omega)$ on the moduli space
of Gieseker stable sheaves defined by a holomorphic
two form $\omega \in H^0(X, K_X)$.  
We will work in the algebraic rather than the K\"ahler category to be able to freely use existing results 
from \cite{HuybrechtsLehn:SheavesBook}. 
See section \ref{SubSection:TauForm} for definitions and notation.  

\begin{Theorem}\label{Theorem:GenericRankTauForm}
Let $X$ be an algebraic surface with $p_g > 0$ and  $b_1 =0$.  
Fix a line bundle $\L \in Pic(X)$ with $c_1(\L) = L$,  and a polarization $H$. 
Let $\omega$ be a holomorphic 2 form.
Suppose the Donaldson series is of the form 
$$
	q_L(X) = \prod_{i=1}^e C_i e^Q (q_0 + q_2 + q_4 + .....)
$$
where the $C_i$ are effective rational divisors, $q_p$ homogeneous polynomials 
of degree $p$ in classes of type $(1,1)$, and $q_0$ a non zero number. Then if $k \gg 0$ 
$$
	\tau(\omega)^n \ne 0 \iff n \le \floor{(d(L,k)-e)/2}
$$
\end{Theorem}

\begin{proof}
First assume that we are in the favorable situation that $H\cdot L$ is odd.
Then semi-stability is the same as stability. Moreover 
there exists a universal sheaf $\F$ over the moduli space 
$\Mbar(\L,k)$ of $H$-stable sheaves with determinant $\L$ and $c_2 = k$
by \cite[Corollary 4.6.7]{HuybrechtsLehn:SheavesBook}. Taking $k$ sufficiently
large, we can assume that the moduli space is generically 
smooth  and reduced of complex dimension $d(L,k)$
\cite[theorem 9.3.3]{HuybrechtsLehn:SheavesBook}. We write $d$ for $d(L,k)$.

The rank of the form $\tau(\omega)$ at the generic point can be determined by finding the maximal dimension
of an even dimensional subvariety on which it is symplectic. In fact we can take the maximum over any family of 
subvarieties that passes through a general point of an open and dense subset of the moduli space. We can therefore
restrict to the varieties $V_H^{(e)}$ cut out by $e$ copies of the nef and big divisor corresponding to 
$\mu(H)$ for $H$ a hyperplane section on $X$ cf. \cite[proof of theorem 2.6]{OGrady:AnaloguesDonaldsonPolynomials} 
\cite[Prop.8.2.3 e.v.]{HuybrechtsLehn:SheavesBook}. With these remarks we see that the non vanishing of $\tau(\omega)^n$ 
is equivalent to the non vanishing of the $L_2$ norm. 
\begin{equation}
\label{Eq:L2norm} 	
	\|\tau(\omega)^n\|_{L_2(V_H^{(d-2n)})}^2 =  
               {\tbinom{2n}{n}}^{-1}\int_{V_H^{(d-2n)}} (\tau(\omega) + \tau(\bar\omega))^{2n}.
\end{equation} 
Here we integrate over the dense smooth locus of $V_H^{(d-2n)}$, 
or any other subset of full measure. Note that this $L_2$ norm depends on the complex structure of $V_H^{(d - 2n)}$, 
(hence of $\Mbar$ and $H$) through the complex conjugation, but is otherwise independent of the choice of a metric. 

\begin{Lemma}\label{Lemma:L2IsDonaldson}
Up to a nonzero constant the $L_2$ norm \eqref{Eq:L2norm} is the value of the Donaldson polynomial
$$
	\|\tau(\omega)^{(d(L,k)-e)/2}\|_{V_H^{(e)}}^2 = \txt{const} q_{L,k}((\omega +\bar \omega)^{d(L,k)-e},H^e).
$$
\end{Lemma}

\begin{proof}
This lemma is essentially the main theorem of O'Grady \cite[theorem 2.4,2.6]{OGrady:AnaloguesDonaldsonPolynomials} 
and Morgan's identification 
of the Donaldson polynomials with their algebraic geometric analogues 
\cite{Morgan:ComparisonAlgebraicDonaldsonInvariants}. Note that 
$d \equiv e \bmod 2$ by assumption.

To be precise, we first have to identify the class $\tau(\omega)$ defined using $\ch_2(\FF)$ with $\mu(\omega)$ 
defined using $1/4 p_1(\FF) = c_2(\FF) - (1/4) c_1^2(\FF)$\footnote
{Alternatively we could define $\tau$ using $p_1(\FF) = ch_2(\sum_i (-1)^i \Ext^i(\FF,\FF))$, where 
$Ext^i(\FF,\FF)$ is the universal $\Ext$ sheaf that exists even when the universal sheaf $\FF$ itself does not. 
However we would then loose contact with \cite{HuybrechtsLehn:SheavesBook}}
.
 
Since these classes differ by a term proportional to $c_1^2(\F)$ we have to show that   
\begin{equation}
\label{Eq:c1TermVanishes}
	\int_X c_1^2(\FF)\wedge \omega = 0.
\end{equation}
To prove this, we decompose $c_1^2(\FF)$ in Hodge type on both $X$ and (the smooth
locus of) $\Mbar$. The only contribution to the ``integral'' over $X$ comes from the 
$(c_1(\FF)^2)^{(02|20)} = (c_1(\FF)^{(01|10)})^2$ part. However 
$c_1(\F)^{(01|10)} = 0$ because we have assumed that $b_1 = 0$\footnote
{to prove \eqref{Eq:c1TermVanishes} It clearly suffices that $H^1(X) \wedge H^1(X) \to H^2(X)$ vanishes (see remark
\ref{Remark:TopologicalCondition})}
. Thus
$$
	\mu(\omega) = \tau(\omega)
$$
as a cohomology class on $\Mbar$.

By an easy extension of 
Morgan's identification of the Donaldson polynomials with their algebraic analogues (up to sign), 
we get writing $d= d(L,k)$, 
\begin{align*}
	q_{L,k}((\omega+\bar\omega)^{(d-e)}, H^e) &= 
         (\mu(\omega) + \mu(\bar\omega))^{d-e} \cup \mu(H)^e\cap [\Mbar ]
\\
         &= (\mu(\omega) + \mu(\bar\omega)^{(d-e)} \cap [V_H^{(e)}]
\\
	   &= (\tau(\omega) + \tau(\bar\omega)^{(d-e)} \cap [V_H^{(e)}]
\\
	   &= \int_{[V_H^{(e)}]}(\tau(\omega) + \bar\tau(\omega))^{(d-e)}
\end{align*}
where the last integral is indeed integration of forms over any smooth open dense locus. This proves the lemma.
\end{proof}

Now to prove theorem \ref{Theorem:GenericRankTauForm}, we plug in the form of the Donaldson series and 
consider the degree $d$ part. 
The cohomology class of $\omega + \bar\omega$ is orthogonal to all $C_i$
and evaluates to 0 on all the homogeneous polynomials $q_p$ in classes of type 
$(1,1)$ if $p\ne 0$.  Therefore the definition of the Donaldson series gives  
$$
	q_{L,k}((\omega + \bar\omega)^{(d-e)}, H^e)  =  
                      q_0 \(\prod_{i=1}^e H\cdot C_i \)\( \int_X \omega\wedge\bar\omega \)^{(d-e)/2}
$$
which is non zero as required. Moreover if $e \ge 2$ we see that 
$$ 
	q_{L,k}((\omega + \bar\omega)^{(d-e) +2 }, H^{(e-2)})  =  0                    
$$
as required. Since the $\tau(\omega)^{(d-e)/2 + 1} = 0$ for trivial dimension reasons if $e <2$, 
the theorem is proved in the case of $H\cdot L$ odd. It remains to lift the restriction $H\cdot L $ odd,
and deal with semi stable sheaves. 

Like in theorem \ref{Theorem:ExistenceStableBundles}, we consider the blow up $\sigma:\~X \to X$ of the surface $X$. 
Fix the class $\~L= L + E$ and a polarization of type 
$\~H = 2\lambda H - E$, where $\lambda \gg 0$ to get back in the $\~H\cdot \~L$ is odd case.
By the blowup formula for the Donaldson series 
(\cite{FintushelStern:BlowupFormulas} or \cite[lemma 9.2, theorem 7.23]{KronheimerMrowka:EmbeddedSurfaces})
we know that 
$$
 	q_{\~L} = \sinh(E)e^{-EE}q_L = E\(\prod_{i=1}^e C_i\) e^{\~Q/2}(q_0 + \~q_2 + \~q_4 + \ldots).
$$
Thus $q_{\~L}$ is divisible by $e+1$ effective divisors, and the $\~q_p$ are homogeneous 
polynomials in classes of type $(1,1)$. 
Then by the above the form $\tau(\sigma^*\omega)^n \ne 0$ for 
$$
	n \le (d(\~L,k) -(e+1))/2 = (d(L,k)-e)/2
$$
Note that  $\tau_{\~X}(\sigma^*\omega)$ has the rank that we claim for $\tau_X(\omega)$, but the 
virtual dimension $d(\~L,k) = d(l,k) +1$ of the moduli space  on which $\tau_{\~X}(\sigma^*\omega)$ lives 
is one larger than for $\tau_X(\omega)$. This has a simple geometric interpretation.

For $k \gg 0$ and $\lambda\gg0$ the moduli space $\Mbar = \Mbar_H(X,\L,k)$ contains a dense open set 
$U$, with  a $\P^1$-bundle $\pi:P\to U$ 
coming with an open immersion $i:P \hookrightarrow \Mbartilde= \Mbar_{\~H}(\~X,\~L,k,)$ 
with dense image by (for the construction see \cite{Brussee:BlownupSurfaces}, 
denseness follows because the moduli space $\Mbartilde$ is irreducible for $k \gg 0$ 
by \cite[Th. 9.4.3]{HuybrechtsLehn:SheavesBook}).
 
Every geometric point in $U$ corresponds to a {\em H-stable} sheaf $\F$ 
without singularity in the blow up point $x$, and the fibre over a point $[\F]$
of $P$ corresponds to the 1 dimensional family of $\~H$ stable extensions
$$
	0 \to \sigma^*\F \to \~\F \to \O_E(E) \to 0 
$$
Fix an extension $\~F_0$ in the family. Since $\O(E)$ has degree $-1$ on the line $E$
it is easy to see that $\sigma_*\~\F_0 \iso \F$, and 
$$
	R_1 \sigma_* \~\F_0 = 0
$$
vanishes.

On $P$ we can pullback $\tau_X(\omega)$ and restrict $\tau_{\~X}(\sigma^*(\omega))$. 
I claim these forms are  equal:
\begin{equation}
\label{Eq:PiTauIsTauSigma}  \pi^* \tau_X(\omega) = \tau_{\~X}(\sigma^*(\omega))|_{P}.
\end{equation}
Since $\tau_{\~X}(\sigma^*\omega)$ has the required rank, and 
$$
	\rank(\tau(X)) = \rank(\sigma^*\tau(X)),
$$
once the claim is proven, we are done. 

To prove \eqref{Eq:PiTauIsTauSigma}, consider the universal family $\~\FF$ on $\Mbartilde$ of 
$\~H$ stable sheaves on $\~X$  with determinant  $\~\L$ and 
second chern class $k$. Its restriction to $P$ is also denoted $\~\F$.
Then the family $(\sigma \times \id)_*\~\FF$ is a family  
of stable sheaves on $X$  parametrised by $P$.  By construction, 
it is constant on every fibre of $P \map{\pi}  U$. Thus the family  
$(\sigma \times \pi)_* \~\FF$ is a family of stable rank 2 sheaves 
on $X$ parametrised by $U$ with determinant $\L$ and $c_2 = k$.
In fact by the interpretation of the $\P^1$ fibre of $\pi$ this family is 
universal on $U$!

Now $\tau_X(\omega)$ is independent of the choice of universal family so, 
\begin{align}
\pi^* \tau(\omega) &=	\pi^* \int_X \ch_2((\sigma\times\pi)_*\~\FF) 
		\wedge \omega 
\\
	&=
    \int_X \ch_2( (\id \times \pi)^*(\sigma\times\pi)_*\~\FF ) \wedge \omega 
\end{align}

We then note that $(\sigma \times \id)_* \~\FF$ is constant along the 
$\P^1$-fibres of $\pi$ to conclude that 
$$
	(\id \times \pi)^*(\sigma \times \pi)_* \~\FF = 
(\id \times \pi)^* (\id \times \pi)_* (\sigma \times \id)_* \~\FF \iso 
        (\sigma \times \id)_*\~\FF
$$
(the most dificult thing here, is to look through the notation). Therefore
$$
	\pi^*\tau(\omega) = \int_X \ch_2\((\sigma \times \id)_* \~\FF\) \wedge \omega.
$$
But all higher derived sheaves $R^i (\sigma \times \id)_* \~\FF = 0$ vanish. 
Therefore by Grothendieck Riemann Roch
$$
	\ch((\sigma \times \id)_* \~\FF) = (\sigma \times\id)_!\(\ch(\~\FF) \td(\sigma\times \id)\)
$$
where $\td(\sigma\times\id) = \td(\~X)/\sigma^*\td(X)$. Thus writing $[\  ]_i$ for taking the $i$ 
dimensional cohomology part we get 
\begin{align*}
	\pi^*\tau_X(\omega) &= \int_X [(\sigma \times \id)_! \(\ch(\~\FF)\td(\~X)\)\wedge \td(X)^{-1}]_4\wedge \omega
\\
		            &= \int_{\~X} [\ch(\~\FF) \wedge \td(\~X) \wedge \sigma^* \td(X)^{-1} \wedge \sigma^*\omega]_6
\\
			      &= \int_{\~X} \ch_2(\~\FF) \wedge \sigma^*\omega
\\
			      &= \tau_{\~X}(\sigma^* \omega)
\end{align*}
where in the third line we used that because the Todd class is a sum of class of type $(p,p)$
the only term in the Todd quotient that survives when multiplied by $\omega$ is the constant term 1. 
The claim is proved.
\end{proof}

\begin{Remark}\label{Remark:TopologicalCondition}
We have used the Kronheimer Mrowka machinery which assumes $b_1 =0$. 
This is probably not essential. What we have used in an essential way, is that the map 
\begin{equation}
\label{Eq:Wedge2Condition}
    \wedge^2 H^1(\O)  \to H^2(\O)
\end{equation}
vanishes, although arguably we could have just used $\mu(\omega)$ to begin with. 
If $b_1 \ne 0$ condition \eqref{Eq:Wedge2Condition} means, that $X$  
admits a fibration to a curve of genus $g \ge 1$ (the map is obtained directly 
from the Albanese map if $b_1 = 2q = 2$, or else by Stein factorization of the map 
defined by a pencil of 1-forms \cite[proposition 4.1]{BarthPetersVanDeVen:ComplexSurfacesBook}.) 

By Hodge theory, 
condition \eqref{Eq:Wedge2Condition} above is equivalent to the purely topological 
condition that the image of 
$$
	\wedge^2 H^1(X,\R) \to H^2(X,\R)
$$
is at most one dimensional. This can be seen as follows. If condition \eqref{Eq:Wedge2Condition} holds then 
the image of $\wedge^2 H^1(X,\R)$ is real, lies in $H^{1,1}$ and is isotropic. Thus by the Hodge index theorem 
the image can be at most one dimensional. On the other hand, if the vanishing condition \eqref{Eq:Wedge2Condition} 
does not hold, the real part of the image of $\wedge^2 H^{01}$ and $\wedge^2 H^{10}$ spans a space 
of dimension at least $\dim_\R \C = 2$. 
\end{Remark}


\end{document}